\documentclass{jT}

\usepackage{amssymb,amsmath,latexsym,amsthm}
\usepackage[french]{babel}
\usepackage[T1]{fontenc}

\newtheorem{thm}{Theorem}[section]
\newtheorem{theoreme}[thm]{Th\'eor\`eme}
\newtheorem{lemme}[thm]{Lemme}

\theoremstyle{definition}

\numberwithin{equation}{section}

\begin{document}

\def\refname{\centerline{Bibliographie}}

\title[Points d'accumulation des nombres de Salem et conjecture de Lehmer]{NOMBRES DE PISOT, NOMBRES DE SALEM ET LA CONJECTURE DE LEHMER}


\author{Mohamed Amara}
\address{Professeur Emérite, Université de Tunis, Tunisie\\
Tunisie}
\email{mehdi.amara@sfr.fr}
\urladdr{http://www.math.u-bordeaux.fr/A2X/} 


\thanks{} 

\maketitle

\begin{resume}
Notre étude consiste à établir le lien entre l'ensemble $S$, des nombres de Pisot, et celui $T$, des nombres de Salem. Le premier lien étroit a été établi par Salem : « tout élément de l'ensemble $S$ est un point d'accumulation de l'ensemble $T$ ». Par le procédé de Boyd, nous montrons que tout point d'accumulation de l'ensemble $T$ est un élément de l'ensemble $S$, ces deux lemmes entraînant que l'ensemble $S \cup T$ est un fermé de la demi-droite réelle $\left]1, +\infty \right[$. On rejette ainsi la conjecture de Lehmer tout en démontrant celle de Boyd.
\end{resume}

\begin{abstr}
We investigate the relationship between the set $S$ of Pisot numbers and the set $T$ of Salem numbers. Salem first established that: " every Pisot number is an accumulation point of the set $T$ ". Building on Boyd’s method, we show that every accumulation point of $T$ belongs to $S$. Together, these results imply that the union $S \cup T$ forms a closed subset of the real half-line $\left]1, +\infty \right[$. Consequently, this settles Boyd’s conjecture while disproving Lehmer’s conjecture.
\end{abstr}

\bigskip

\section{Introduction}
Dans cet article, nous proposons d'étudier les liens entre deux ensembles d'entiers algébriques, celui $S$ des nombres de Pisot et celui $T$ des nombres de Salem. Ces deux ensembles jouent un rôle fondamental en théorie des nombres et, reprenant l'introduction de l'article de Dufresnoy et Pisot \cite{Dufresnoy1953}, permettent « d'obtenir des relations nouvelles entre corps différents [...] et sont en liaison étroite avec l'approximation rationnelle, donc se rattachent aussi à certaines parties de la théorie des séries trigonométriques ». Le premier lien étroit entre ces deux ensembles revient, par un procédé appartenant à Salem \cite{Salem1944}, à établir que « tout élément de l'ensemble de $S$ de Pisot est un point d'accumulation de l'ensemble $T$ de Salem ». D'après Lehmer \cite{Lehmer1933}, la réciproque serait impossible.\newline
Nous avons précédemment étudié le sous-ensemble dérivé des nombres de Pisot $S^\prime(1)$ des points d'accumulation de $S$ entre 1 et 2 \cite{Amara1966}. On en propose ici une généralisation à partir du sous-ensemble dérivé $S^\prime(u_0)$, où $u_0$ est un entier supérieur ou égal à 2.  Par un procédé de Boyd \cite{Boyd1977}, réciproque de celui de Salem, nous montrerons que « tout point d'accumulation de l'ensemble $T$ appartient à l'ensemble $S$ », ce qui permet de conclure que l'ensemble $S \cup T$ est un fermé de la demi-droite $\left]1, +\infty \right[$, de plus petit $\tau_0$, zéro du polynôme \[1+z-z^3\left(1+z+z^2+z^3+z^4\right)+z^9+z^{10}\] Ce résultat permet de rejeter la conjecture de Lehmer \cite{Lehmer1933} et confirme celle de Boyd \cite{Boyd1977}. 

\section{Sur un sous-ensemble dérivé des nombres de Pisot}

\subsection{Définitions}\hfill\\
On entend par nombre de Pisot un entier algébrique $\theta$ strictement supérieur à 1, dont tous les conjugués sont situés dans la région $\left|z\right|<1$. Le polynôme, à coefficients entiers, dont l'élément $\theta$ est zéro, est désigné par la lettre $P$ s’écrivant
\begin{equation*}
P(z)=p_0+p_1z+\ldots+p_{s-1}z^{s-1}+\varepsilon z^s \text{ avec }p_0\geq1 , \varepsilon=\pm1
\end{equation*}
Le polynôme $P$ est irréductible puisque, s’il se décomposait en facteurs, l'un des facteurs aurait tous ses zéros dans la région $\left|z\right|<1$ alors que le terme de plus haut degré est de module égal à 1, ce qui est impossible. L’ensemble des nombres de Pisot est désigné par la lettre $S$. Rappelons que cet ensemble est un fermé de la droite réelle \cite{Dufresnoy1953}. Le plus petit élément de l'ensemble $S$ est le zéro, noté $\theta_0$, du polynôme 
\begin{equation*}
P (z) = 1 + z-z^3
\end{equation*}
On désigne par la lettre $S'$ l'ensemble des éléments d’accumulation de l'ensemble $S$. Les éléments de $S'$ sont caractérisés dans \cite{Dufresnoy1953} par :
\begin{quotation}
« Pour qu'un nombre de l'ensemble $S$ appartienne à l’ensemble $S’$, il faut et il suffit qu'il existe au moins un polynôme $A$, à coefficients entiers tel qu'on ait sur $\left|z\right|=1$, $\left|A(z)\right|\le\left|P(z)\right|$, l’égalité n'ayant lieu qu'en un nombre fini de points. »
\end{quotation}
Le plus petit élément de l'ensemble $S’$ est le nombre $\theta'$ zéro du polynôme
\begin{equation*}
P(z)=1+z-z^2
\end{equation*}
associé aux polynômes
\begin{equation*}
A(z)=1\ ;A(z)=1-z^2
\end{equation*}

Aux polynômes $A$ et $P$ représentant un élément de l’ensemble $S$, on adjoint les polynômes suivants :
\begin{multline*}
B\left( z \right)=\varepsilon^\prime z^h A\left(1/z\right) , \text{\textit{h}= degré de }A, \varepsilon'=\pm1 \text{ choisi tel que } B(0)=b_0 \geqslant 1\\
Q(z)=\varepsilon z^s P\left( 1/z \right) ; \varepsilon=\pm1\ \text{ choisi tel que } Q(0)=1
\end{multline*}

\subsubsection{Caractérisation des éléments de $S^\prime$}\hfill\\
La caractérisation des éléments de l’ensemble $S^\prime$ est illustrée par le polynôme $G$ ainsi défini :
\paragraph{1/ Pour $h\geq s+1$ }
\[G(z)=A(z)B(z)-\varepsilon\varepsilon^\prime z^{h-s}P(z)Q(z)\]
	tel que sur $\left|z\right|=1$
\[G(z)=-\varepsilon^\prime z^h\left[\left|Q(z)\right|^2-\left|A(z)\right|^2\right]\]

\paragraph{2/ Pour $s\geq h+1$}
\[G(z)=P(z)Q(z)-\varepsilon\varepsilon^{\prime} z^{s-h}A(z)B(z)\]
	tel que sur $\left|z\right|=1 $
\[G(z)=\varepsilon z^s\left[\left|Q(z)\right|^2-\left|A(z)\right|^2\right]\]

\paragraph{3/ Pour $s=h$}
\[\eta z^kG(z)=A(z)B(z)-\varepsilon\varepsilon^\prime P(z)Q(z)\]
tel que sur $\left|z\right|=1$
\[\eta z^kG(z)=-\varepsilon^\prime z^s\left[\left|Q(z)\right|^2-\left|A(z)\right|^2\right]\]

L’entier $k\geq0$ est choisi tel que $G(0)\neq0$ ; le coefficient $\eta=\pm1$ est choisi tel que
$G(0)\geq1$.\newline

L’inégalité sur $\left|z\right|=1$
\begin{equation*}
\left|A(z)\right|\le\left|Q(z)\right|
\end{equation*}
n’ayant lieu qu’en un nombre fini de points, entraine deux propriétés pour le polynôme $G$ qui seront systématiquement utilisées tout le long de notre étude à savoir :
\begin{enumerate}
\item[i)] Le polynôme $G$ ne peut-être identiquement nul
\item[ii)] Sur $\left|z\right|=1$ le polynome $G$ quand il n’est pas nul admet un argument du type $\varepsilon'' z^l$. Il en résulte que tout zéro de $G$, sur $\left|z\right|=1$, est de multiplicité paire.
\end{enumerate}
Pour tout entier $u_0\geq2$, on désigne par :
\begin{equation*}
S^\prime(u_0)=S^\prime \cap\left]u_0,u_0+1\right[
\end{equation*}
un sous-ensemble dérivé de $S’$ tel que les polynômes $A$ et $P$ représentants tout élément de $S^\prime(u_0)$ vérifient
\begin{equation*}
A(0)=P(0)=u_0
\end{equation*}

\subsection{Plus petits éléments de $S'(u_0)$}
\subsubsection{Étude du cas $h\geq s+1$}\hfill\\
A partir des différents polynômes associés à l’élément $\theta$ de $S^\prime (u_0)$ formons les polynômes :
\begin{eqnarray*}
A_1(z)&=&A^2(z)+G(z)=A^2(z)+A(z)B(z)-\varepsilon\varepsilon^\prime z^{h-s}P(z)Q(z)\\
B_1(z)&=&B^2(z)+G(z)=B^2(z)+A(z)B(z)-\varepsilon\varepsilon^\prime z^{h-s}P(z)Q(z)\\
G_1(z)&=&A_1(z)B_1(z)-z^{2(h-s)}\left[P(z)Q(z)\right]^2
\end{eqnarray*}

Les polynômes $A_1$, $B_1$ vérifient sur $\left|z\right|=1$
\begin{equation*}
\left|A_1(z)\right|=\left|B_1(z)\right|\le\left|A(z)\right|^2+\left|G(z)\right|^2=\left|Q(z)\right|^2
\end{equation*}
L’égalité $\left|A_1(z)\right|=\left|Q(z)\right|^2$ n’a lieu qu’en un nombre fini de points sur $\left|z\right|=1$, autrement le polynôme $G$ serait identiquement nul. De plus, l’inégalité $\left|A_1(z)\right|\le\left|Q(z)\right|^2$ sur $\left|z\right|=1$,  justifie la formation du polynôme $G_1$. Le procédé par lequel nous avons introduit les polynômes $A_1$, $B_1$ et $G_1$, est reproduit dans la récurrence 
\begin{eqnarray*}%
A_{n+1}(z)&=&A_n^2(z)+G_n(z)=A_n^2(z)+A_n(z)B_n(z)-z^{2^n(h-s)}\left[P(z)Q(z)\right]^{2^n}\\
B_{n+1}(z)&=&B_n^2(z)+G_n(z)=B_n^2(z)+A_n(z)B_n(z)-z^{2^n(h-s)}\left[P(z)Q(z)\right]^{2^n}\\
G_{n+1}(z)&=&A_{n+1}(z)B_{n+1}(z)-z^{2^{n+1}\left(h-s\right)}\left[P(z)Q(z)\right]^{2^{n+1}}\\&=&\left[A_n(z)B_n(z)\right]^2G_n(z)
\end{eqnarray*}

Les polynômes $A_{n+1}$, $B_{n+1}$ vérifient sur $\left|z\right|=1$
\begin{equation*}
\left|A_{n+1}(z)\right|=\left|B_{n+1}(z)\right|\le\left|A_n(z)\right|^2+\left|G_n(z)\right|^2=\left|Q(z)\right|^{2^{n+1}}
\end{equation*}
L’égalité $\left|A_{n+1}(z)\right|=\left|Q(z)\right|^{2^{n+1}}$ sur $\left|z\right|=1$ n’a lieu qu’en un nombre fini de points, autrement le polynôme $G_{n+1}$ serait identiquement nul, ce qui n’est pas le cas.\newline
Considérons la suite de fonctions
\begin{equation*}
f_n(z)=\frac{A_n(z)+B_n(z)}{\left[P(z)Q(z)\right]^{2^{n-1}}}
\end{equation*}
Ces fonctions $(f_n)$ sont non constantes, méromorphes dans $\left|z\right|<1$, de module sur  $\left|z\right|=1$, inférieur ou égal à 1. Les relations de récurrence vérifiées par les polynômes $A_n$ et $B_n$ entraînent pour les suites $(f_n)$ et $(a_n=f_n(0))$ les récurrences :
\begin{equation*}
f_{n+1}(z)=\left[f_n(z)\right]^2-2z^{2^n\left(h-s\right)}\ ;\ a_{n+1}=a_n^2\ ,\ n\geq1
\end{equation*}
Formons la suite de fonctions
\begin{equation*}
F_n(z)=f_n(z)\left[\frac{P(z)}{Q(z)}\right]^{2^{n-1}}\left(\frac{1-\theta z}{\theta-z}\right)^{2^n}=\frac{A_n(z)+B_n(z)}{\left[Q(z)\right]^{2^n}}\left(\frac{1-\theta z}{\theta-z}\right)^{2^n}
\end{equation*}
La fonction $F_n$ est holomorphe dans $\left|z\right|<1$, de module, sur $\left|z\right|=1$, inférieur ou égal à 2. Par le lemme de Schwartz
\begin{equation*}
F_n(0)=f_n(0)\left[\frac{u_0}{\theta^2}\right]^{2^{n-1}}=a_n\left[\frac{u_0}{\theta^2}\right]^{2^{n-1}}\le2
\end{equation*}
La suite $(a_n)$ s’écrivant
\begin{equation*}
a_n=\left(a_1\right)^{2^{n-1}}\le2
\end{equation*}
On obtient :
\begin{equation*}
F_n(0)=\left[\frac{a_1u_0}{\theta^2}\right]^{2^{n-1}}\le2
\end{equation*}
Soit
\begin{equation*}
\frac{a_1u_0}{\theta^2\ }\le{(2)}^\frac{1}{2^{n-1}}
\end{equation*}
Par passage à la limite avec $n$
\begin{equation}
\label{equ1}
\theta^2\geq a_1u_0=u_0f_n(0)=A_1(0)+B_1(0)=u_0^2+u_0+b_0^2+u_0\geq\left(u_0+1\right)^2
\end{equation}
Soit
\begin{equation*}
\theta\geq u_0+1
\end{equation*}
Or l’élément $\theta$ appartient à $S^\prime (u_0)$ et vérifie
\begin{equation*}
\theta<u_0+1
\end{equation*}
L’inégalité $h\geq s+1$ est à rejeter 

\subsubsection{Etude du cas $s\geq h$}\hfill\\
Notons tout d’abord
\paragraph{1/ Pour $s\geq h+1$}
\begin{equation*}
G(z)=z^{2s}G\left(\frac{1}{z}\right);G(0)\geq u_0
\end{equation*}
\paragraph{2/ Pour $s=h$}
	\begin{itemize}
	\item[$\bullet$] L’entier $k=0$ entraine
		\begin{equation*}
		G(z)=z^{2s}G\left(\frac{1}{z}\right);G(0)=u_0\left(b_0-1\right)\geq1, \text{ soit } b_0\geq2\ \text{et}\ G(0)\geq u_0
		\end{equation*}
	\item[$\bullet$] L’entier $k\geq1$ entraine
		\begin{equation*}
		G(z)=z^{2\left(s-k\right)}G\left(\frac{1}{z}\right)
		\end{equation*}
	\end{itemize}

Introduisons comme précédemment les polynômes 
\begin{eqnarray*}
&A_1(z)=&A^2(z)+G(z)\\
&B_1(z)=&z^{2s}A_1\left(\frac{1}{z}\right)=\left\{\begin{array}{c}G(z)+z^{2\left(s-h\right)}B^2(z),\ s\geq h+1 \\B^2(z)+z^{2\left(s-k\right)}G\left(\frac{1}{z}\right),\ s=h\ \text{et}\ k\geq0\end{array}\right.\nonumber \\
&G_1(z)=&A_1(z)B_1(z)\left[P(z)Q(z)\right]^2\nonumber\\
&=&\left\{\begin{array}{c}\left[A(z)-\varepsilon\varepsilon^\prime z^{s-h}B(z)\right]^2G(z),\ s\geq h+1 \\\left[B(z)-\eta z^kA(z)\right]^2G(z),\ s=h\end{array}\right.\nonumber
\end{eqnarray*}
De par son expression, le polynôme $G_1$ ne peut être identiquement nul. De plus les polynômes $A_1$ et $B_1$ vérifient sur $\left|z\right|=1$
\begin{equation*}
\left|A_1(z)\right|=\left|B_1(z)\right|\le\left|A(z)\right|^2+\left|G(z)\right|=\left|Q(z)\right|^2
\end{equation*}
L’égalité $\left|A_1(z)\right|= \left|Q(z)\right|^2$ sur $\left|z\right|=1$ n’a lieu qu’en un nombre fini de points. Le procédé d’itération utilisé dans le cas $h\geq s+1$ permet d’obtenir des suites de polynômes ${(A}_n),\ \left(B_n\right)$ et ${(G}_n)$ vérifiant
\begin{eqnarray*}
A_{n+1}(z)&=&A_n^2(z)+A_n(z)B_n(z)-\left[P(z)Q(z)\right]^{2^n}\\
B_{n+1}(z)&=&B_n^2(z)+A_n(z)B_n(z)-\left[P(z)Q(z)\right]^{2^n}\\
G_{n+1}(z)&=&\left[A_n(z)+B_n(z)\right]^2G_n(z)
\end{eqnarray*}
Il résulte de la récurrence que le polynôme $G_{n+1}$ ne peut pas être identiquement nul. Les polynômes $A_{n+1}$, $B_{n+1}$ vérifient sur $\left|z\right|=1$
\begin{equation*}
\left|A_{n+1}(z)\right|=\left|B_{n+1}(z)\right|\le\left|A_n(z)\right|^2+\left|G_n(z)\right|=\left[Q(z)\right]^{2^{n+1}}
\end{equation*}
L’égalité $\left|A_{n+1}(z)\right|=\left[Q(z)\right]^{2^{n+1}}$ sur $\left|z\right|=1$ n’a lieu qu’en un nombre fini de points.\newline
Considérons la suite de fonctions rationnelles 
\begin{equation*}
f_n(z)=\frac{A_n(z)+B_n(z)}{\left[P(z)Q(z)\right]^{2^{n-1}}}\ ,\ n\geq1
\end{equation*}
Ces fonctions sont non constantes, méromorphes dans la région $\left|z\right|<1$ de module sur $\left|z\right|=1$, inférieur ou égal à 2.\newline
Les suites $\left(f_n(z)\right)$, $(a_n=f_n(0))$ vérifient les récurrences 
\begin{equation}
\label{equ2}
f_{n+1}=f_n^2(z)-2,\ a_{n+1}=a_n^2-2,\ n\geq1
\end{equation}
Étudions la suite $(a_n)$ de l'équation (\ref{equ2})
\begin{equation*}
u_1 a_1 = A_1(0)+B_1(0)=\left\{\begin{array}{l} u_0^2 + 2u_0 \text{ pour } s\geq h+1 \\u_0^2+b_0^2+u_0 \left( b_0 -1 \right) \text{ pour } k=0, s+h, b_0\geq2 \\ u_0^2+b_0^2+G(0) \text{pour} k\geq1, s=h, b_0\geq1\end{array}\right.
\end{equation*}
On obtient dans tous les cas
\begin{equation*}
u_0 a_1\geq u_0^2+2 \text{ ou } a_1>u_0+1>2
\end{equation*}
La suite $(a_n)$, strictement croissante vérifie pour tout $n$
\begin{equation*}
a_n>2
\end{equation*}
Il existe une suite $\left(c_n\right)$ de nombres positifs telle que
\begin{equation*}
a_n=2\;ch2^n c_n,\ n\geq1
\end{equation*}
Des propriétés de la fonction hyperbolique et de la récurrence formée par la suite $(a_n)$, on établit
\begin{equation*}
ch2^{n+1} c_{n+1}=2\;ch(2^nc_n)-1=ch2^{n+1}\;c_n
\end{equation*}
La suite $(c_n)$ est stationnaire de valeur commune déterminée par l’égalité
\begin{equation}
\label{equ3}
a_1=2\;ch2c \text{ ou } e^{4c}-a_1 e^{2c}+1=0
\end{equation}

Reprenons la fonction
\begin{equation*}
F_n(z)=f_n(z)\left[\frac{P(z)}{Q(z)}\right]^{2^{n-1}}\left(\frac{1-\theta z}{\theta-z}\right)^{2^n}\ ,\ n\geq1
\end{equation*}
De lemme de Schwarz 
\begin{equation*}
F_n(0)=a_n\left(\frac{u_0}{\theta^2}\right)^{2^{n-1}}\le2
\end{equation*}

Soit pour $a_n$ l’encadrement
\begin{equation*}
e^{2^nc}\le a_n\le2\left[\frac{\theta^2}{u_0}\right]^{2^{n-1}}\le2
\end{equation*}
Par passage à la limite avec $n$ on obtient
\begin{equation*}
\theta^2\geq u_0e^{2c}
\end{equation*}
En revenant à l'équation \ref{equ3} et à $s\geq h$
\begin{equation}
\label{equ4}
\theta^4-u_0a_1\theta^2+u_0^2\geq0
\end{equation}
\newline

\paragraph{1/ Etude du cas $s\geq h+1$}\hfill\\
De l'équation (\ref{equ2}) avec $G(0)=u_0$, on tire 
\begin{equation*}
u_0a_1=u_0^2+2u_0
\end{equation*}
Et de l'équation (\ref{equ4}) on obtient
\begin{equation*}
\theta^4-(u_0^2+2u_0{)\theta}^2+u_0^2\geq0
\end{equation*}
Ou plus simplement
\begin{equation*}
\theta^2-u_0\theta-u_0\geq0
\end{equation*}
Le plus petit élément de $S^\prime(u_0)$, pour $s\geq h+1$, est le zéro du polynôme
\begin{equation*}
P(z)=u_0+u_0z-z^2
\end{equation*}
Le polynôme $A$ n’est autre que
\begin{equation*}
A(z)=u_0+\left(u_0-1\right)z
\end{equation*}
Pour la suite de notre étude, le plus petit élément de $S^\prime(u_0)$, pour $s=h+1$ sera noté $\beta_2$ et représenté par la fonction de rang infini $A(z)/Q(z)$ où 
\begin{eqnarray}
\label{equ5}
\left(1-z\right)A(z)&=&u_0-z-\left(u_0-1\right)z^2\\
(1-z)Q(z)&=&1-\left(u_0+1\right)z+u_0z^3\nonumber\\
G(z)&=&u_0\left(1-z^2\right)^2\nonumber
\end{eqnarray}
\newline

\paragraph{2/ Etude du cas $s=h$}\hfill\\
Pour $k=0$ on dispose de
\begin{equation*}
b_0\geq2,\ G(z)=z^{2s}G\left(\frac{1}{z}\right),\ G(0)=u_0(b_0-1)\geq u_0
\end{equation*}
De l'équation (\ref{equ1}) on tire
\begin{equation*}
u_1a_1\geq u_0^2+2u_0+H=\left(u_0+2\right)^2-2u_0
\end{equation*}
De l'équation (\ref{equ4}) l’élément $\theta$ vérifie
\begin{equation*}
\theta^4-\left[\left(u_0+2\right)^2-2u_0\right]\theta^2+u_0^2\geq0
\end{equation*}
Ou plus simplement
\begin{equation*}
\theta^2-\left(u_0+2\right)\theta+u_0\geq0
\end{equation*}

Et par suite
\begin{equation*}
\theta>u_0+1
\end{equation*}
La valeur $k=0$ est à rejeter.\newline
Pour $k\geq1$ on dispose
\begin{equation}
\label{equ6}
G(z)=z^{2\left(s-k\right)}G\left(\frac{1}{z}\right),\ G(0)\geq1,\ b_0=1,\ \varepsilon\varepsilon^\prime=1
\end{equation}
De l'équation (\ref{equ1}) on tire
\begin{equation*}
u_0a_1=u_0^2+G(0)+1
\end{equation*}
Et de l'équation (\ref{equ4}) l’élément $\theta$ vérifie
\begin{equation}
\label{equ7}
\theta^4-\left(u_0^2+G(0)+1\right)\theta^2+u_0^2>0
\end{equation}
Notons que pour $u_0=1$ et $G(0)=1$ on retrouve l’élément $\theta^\prime$, plus petit élément de $S’$ et zéro du polynôme
\begin{equation*}
1+z-z^2
\end{equation*}
Évaluons le nombre $G(0)$ à partir des développements à l’origine et à coefficients des entiers des fonctions, représentants l’élément $\theta$ de $S^\prime(u_0)$
\begin{equation*}
\frac{A(z)}{Q(z)}=u_0+u_1z+\ldots;\ \frac{B(z)}{Q(z)}=1+b_1z+\ldots;\ \frac{P(z)}{Q(z)}=u_0+v_1z+\ldots
\end{equation*}
On obtient
\begin{equation}
\label{equ8}
G(0)=u_1-v_1+b_1u_0 \text{ avec } b_1\geq1
\end{equation}
Reprenons les écritures des polynômes $A_1$, $B_1$ avec $k=1$
\begin{eqnarray*}
A_1(z)&=&A^2(z)+G(z)=\left[u_0^2+G(0)\right]+\ldots\\
B_1(z)&=&B^2(z)+z^2G(z)=1+\ldots+G(0)z^2
\end{eqnarray*}
Pour déterminer le plus petit élément de $S^\prime(u_0)$ pour $s=h$, nous devons disposer de l’écriture de $A_1$ jusqu’à l’ordre 2 ce qui semble se heurter à de sérieuses difficultés. Toutefois nous retenons que la détermination de plus petit élément de $S^\prime(u_0)$ dépend de la plus petite valeur de $G(0)$.\newline

\subsection{Sur la composition de l’ensemble $S’(u_0)$} \label{I-2}

\subsubsection{Étude du cas $s\geq h+1$}\hfill\\
De la représentation du plus petit élément $\beta_2$ de l’ensemble $S^\prime(u_0)$ et de celle de la suite $S^\prime(1)$, voir Chapitre III de \cite{Amara1966}, à savoir
\begin{equation*}
\left(1-z\right)A(z)=1-z^{s-1};\left(1-z\right)Q(z)=1-2z+z^{s+1},\ s\geq2
\end{equation*}
et par un procédé heuristique on établit
\begin{lemme}\label{lemme1}
Pour $u_0\geq2$ et $ s\geq h+1$ seule la suite $(\beta_s)$ représentée par la famille de fonctions de rang infini où
\begin{equation*}
\left(1-z\right)A(z)=u_0-z^{s-1}-\left(u_0-1\right)z^s;\left(1-z\right)Q(z)=1-\left(u_0+1\right)z+u_0z^{s+1},\ s\geq2
\end{equation*}
appartient à l’ensemble $S^\prime(u_0)$. De plus cette famille vérifie 

\begin{eqnarray*}
Q(z)+zA(z)=P(z)-B(z)=1-z^s\\
G(z)=u_0\left(1-z^s\right)^2
\end{eqnarray*}
L’expression du polynôme $G$ provient du système
\begin{eqnarray*}
Q(z)+zA(z)=1-z^s\\
P(z)-B(z)=1-z^s
\end{eqnarray*}
duquel on obtient 
\begin{equation*}
G(z)=\left(1-z^s\right)\left(B(z)+Q(z)\right)=u_0\left(1-z^s\right)^2
\end{equation*}
\end{lemme}

\subsubsection{Étude du cas $s=h$}\hfill\\
Dans cette partie de notre étude, l’intérêt sera porté en premier lieu vers le développement de la fonction de rang infini
\begin{equation*}
\frac{B(z)}{Q(z)}=1+b_1z+b_2z^2+\ldots
\end{equation*}
Rappelons de l'équation (\ref{equ6})
\begin{equation*}
\eta=1,\ k=1,\ \varepsilon\varepsilon^\prime=1\ et\ b_1\geq1
\end{equation*}
Avec ces données, il existe un polynôme $V$ à coefficients entiers tel que 
\begin{eqnarray}
\label{equ9}
B(z)-Q(z)&=&zV(z)\\
A(z)-P(z)&=&\varepsilon\varepsilon_1z^{s-d-1}V_1(z)\nonumber
\end{eqnarray}
Avec $V(0)=b_1$, $V_1(z)=\varepsilon_1 z^d V(\frac{1}{z})$, $d$ = degré de $V \leq s-2$, $\varepsilon_1 = \pm 1$ choisi tel que $V(0)\geq1$
Par le théorème de Rouché, le polynôme $V$ admet tous ses zéros dans la région $\left|z\right|\geq1$. Par le lemme de Schwarz
\begin{equation*}
1\le V_1(0)\le V(0)=b_1
\end{equation*}
Examinons la relation 
\begin{equation*}
zG(z)=A(z)B(z)-P(z)Q(z)=z\ \left[A(z)V(z)+\varepsilon\varepsilon_1z^{s-d-2}Q(z)V_1(z)\right]
\end{equation*}
\begin{itemize}
\item[$\bullet$] Pour $d=s-2$
\begin{equation*}
G(0)=u_0 b_1+\varepsilon\varepsilon_1 V_1 (0)
\end{equation*}
\item[$\bullet$] Pour $d\le s-3$
\begin{equation*}
G(0)= u_0 V(0) = u_0 b_1
\end{equation*}
\end{itemize}
\hfill\\
Rappelons du commentaire de l'équation (\ref{equ6}) que le plus petit élément de l’ensemble $S^\prime(u_0)$, pour $s=h$, dépend de la plus petite valeur de $G(0)$ soit
\begin{itemize}
\item[$\bullet$] Pour $d=s-2$
\begin{equation*}
G(0)=u_0 -1
\end{equation*}
\qquad avec $b_1=1$, $ \varepsilon\varepsilon_1=-1$, $V_1(0)=1$ et $ u_1-v_1=-1$
\item[$\bullet$] Pour $d\le s-3$
\begin{equation*}
G(0)=u_0
\end{equation*}
\qquad avec $b_1=1$ , $u_1-v_1=0$\newline
\qquad Avec la valeur $b_1=1$ on obtient
\begin{equation}
\label{equ10}
\frac{B(z)}{Q(z)}=1+z+b_2z^2+\ldots
\end{equation}

\end{itemize}
A l’instar de ce développement examinons celui de la fonction
\begin{equation*}
\frac{A(z)}{Q(z)}=u_0+u_1z+\ldots
\end{equation*}
D'après Dufresnoy-Pisot \cite{Dufresnoy1955} associons à ce développement les polynômes
\begin{equation*}
D_2(z)=u_0+\frac{u_1}{u_0+1}z-z^2;D_2^\ast(z)=u_0-\frac{u_1}{u_0-1}z+z^2
\end{equation*}
Possédant chacun un unique zéro $\tau_2$, $\tau_2^\ast$ dans $\left|z\right|>1$ encadrant l’élément $\theta$
\begin{equation*}
1<\tau_2<\theta<\tau_2^\ast
\end{equation*}
\newline
De l’encadrement
\begin{equation*}
\tau_2\le\theta<u_0+1
\end{equation*}
On obtient
\begin{equation*}
D_2\left(u_0+1\right)=u_1-(u_0^2+u_0+1)<0
\end{equation*}
\qquad Soit, pour l’entier $u_1$
\begin{equation*}
u_1\le u_0^2+u_0
\end{equation*}
De l’inégalité 
\begin{equation*}
\theta<\tau_2^\ast
\end{equation*}
On tire 
\begin{equation*}
\frac{1}{\theta}D_2^\ast\left(\theta\right)=\frac{u_0}{\theta}+\theta-\frac{u_1}{u_0-1}<0
\end{equation*}
Soit
\begin{equation*}
u_1>(u_0-1)(\theta+\frac{u_0}{\theta})
\end{equation*}
Ou
\begin{equation*}
u_1^2>\left(u_0-1\right)^2\left[\left(\theta-\frac{u_0}{\theta}\right)^2+4u_0\right]
\end{equation*}
Revenons à l'équation (\ref{equ7}) avec $G(0)=u_0-1$
\begin{equation*}
\theta^4-\left(u_0^2+2u_0-u_0\right)\theta^2+u_0^2
\end{equation*}
Où
\begin{equation*}
(\theta-\frac{u_0}{\theta})\geq u_0(u_0-1)
\end{equation*}
Et, par suite
\begin{equation*}
u_1^2\geq\left(u_0-1\right)^2\left[u_0^2+3u_0\right]
\end{equation*}
Ou mieux
\begin{equation*}
u_1^2\geq\left(u_0-1\right)^2(u_0+2)\geq\left(u_0-1\right)^2(u_0^2+3u_0)
\end{equation*}
On conclut pour le coefficient $u_1$ :
\begin{align}
\label{equ11}
\text{– Pour } G(0)&=u_0-1\text{, on a } u_1=u_0^2+u_0-2 \\
\text{– Pour } G(0)&=u_0\text{, on a } u_0^2+u_0-1 \le u_1 \le u_0^2+u_0 \nonumber  
\end{align}

Des développements des fonctions  $A(z)/Q(z)$, $B(z)/Q(z)$ on établit
\begin{lemme}
\label{lemme3}
\hfill
\begin{itemize}
\item[1/] Pour $u_0\geq2$,  $d\le s-2$
\begin{equation*}
u_1=u_0^2+u_0-2\ ,\ v_1=u_0^2+u_0-1\ ,\ b_1=1\ ,\ b_2=u_0+1\ ,\ G(0)=u_0
\end{equation*}
\item[2/] Pour $u_0\geq2$, $0\le d\le s-3$
\begin{equation*}
u_0^2+u_0-1\le u_1\le u_0^2+u_0\ ,\ u_1=v_1\ ,\ b_1=1\ ,\ b_2=u_0+1\ ,\ G(0)=u_0
\end{equation*}
\end{itemize}
\end{lemme}

Avec les données de 1/ et 2/ la fonction
\begin{equation*}
f(z)=\frac{B(z)+zP(z)}{Q(z)+zA(z)}=1+z+z^2+\ldots
\end{equation*}
admet le même développement tant pour $d=s+2$ que pour $d\le s-3$.
Formons la fonction
\begin{equation*}
f(z)=\frac{B(z)+zP(z)}{Q(z)+zA(z)}=\frac{B(z)}{Q(z)}-z^2\frac{G(z)}{Q(z)(Q(z)+zA(z))}
\end{equation*}
Du développement de l'équation (\ref{equ10}) on obtient celui de $f$
\begin{equation}
\label{equ12}
f(z)=1+z+\left[b_2-G(0)\right]z^2+\ldots
\end{equation}
auquel on associe le polynôme
\begin{equation*}
D_2(z)=1+\frac{1}{2}z-z^2 \text{, avec } W_2=\frac{1}{2}, \; b_2-G(0)\geq\frac{1}{2}
\end{equation*}
soit, pour l’entier $b_2$, la minoration
\begin{equation*}
b_2\geq G(0)+1
\end{equation*}

Au développement de l'équation (\ref{equ10}) associons le polynôme
\begin{equation*}
D_3(z)=1-\left(b_2-2\right)z+\left(b_2-1\right)z^2-z^3
\end{equation*}
Vérifiant
\begin{equation*}
D_3(u_0+1)<0
\end{equation*}
D’où on déduit
\begin{equation*}
b_2\le u_0+2
\end{equation*}
L’entier $b_2$ vérifie l’encadrement
\begin{equation*}
G(0)+1\le b_2\le u_0+2
\end{equation*}

\subsubsection{Etude du cas $d=s-2$, $G(0)=u_0-1$} \label{equ13}\hfill\\
Des valeurs $V(0)=b_1=1$, $ \varepsilon\varepsilon_1=-1$ , $V_1(0)=1$ , $ u_1=v_1=-1$,
 le système (\ref{equ9}) s’écrit
\begin{eqnarray}
\label{equ14}
B(z)-Q(z)=zV(z)\\
A(z)-P(z)=-zV(z) \nonumber
\end{eqnarray}
De Dufresnoy-Pisot (Lemme 3 dans \cite{Dufresnoy1953}), les zéros du polynôme $V$ sont tous situés
sur $\left| z\right|=1$ et tous distincts.
De la multiplicité des zéros de $G$ sur $\left|z\right|=1$, il existe un polynôme $K$ tel que
\begin{equation}
\label{equ15}
G(z)=\left[P(z)-B(z)\right]V(z)=K(z)\left[V(z)\right]^2
\end{equation}

Du degré de $V=s-2$, de celui de $G=2(s-1)$, le polynôme $K$ s’écrit
\begin{equation*}
K(z)=\left(u_0-1\right)+kz+\left(u_0-1\right)z^2
\end{equation*}

Des développements dans l’identité
\begin{equation*}
\frac{P(z)}{Q(z)}-\frac{B(z)}{Q(z)}=K(z)\left[\frac{V(z)}{Q(z)}\right]
\end{equation*}

On obtient
\begin{equation*}
u_0^2+u_0-2=v_1-1=k+\left(u_0-1\right)b_2
\end{equation*}
L’entier $k$ est un multiple de $\left(u_0-1\right)$, soit
\begin{equation*}
k=\left(u_0-1\right)k_1\ ,\ u_0+2=k_1+b_2
\end{equation*}

Les zéros du polynôme $K$, sur $\left|z\right|=1$, sont de multiplicité paire comme le sont ceux du polynôme $G$ soit
\begin{equation*}
k_1=2\ ,\ k=2\left(u_0-1\right)\ ,\ b_2=u_0
\end{equation*}
et on conclut
\begin{equation*}
f(z)=1+z+z^2+\ldots
\end{equation*}
et de l'équation (\ref{equ15}) on obtient
\begin{equation}
\label{equ16}
A(z)-Q(z)=P(z)-B(z)=\left(u_0-1\right)\left(1+z\right)^2V(z)
\end{equation}

\subsubsection{Etude du cas $d\le s-3$,  $G(0)=u_0$}\hfill\\
A partir du paragraphe \ref{equ13}, le coefficient $b_2$ vérifie l’encadrement
\begin{equation*}
u_0+1\le b_2\le u_0+2
\end{equation*}

Examinons le développement
\begin{equation*}
\frac{B(z)}{Q(z)}=1+z+\left(u_0+2\right)z^2+\ldots
\end{equation*}

Auquel est associé le polynôme
\begin{equation*}
D_3(z)=1-u_0z+\left(u_0+1\right)z^2-z^3\ ;w_3=u_0^2+2u_0+2
\end{equation*}
Du développement de la fonction
\begin{eqnarray*}
\frac{B(z)+z^2P(z)}{Q(z)+z^2A(z)}&=&\frac{B(z)}{Q(z)}-z^3\frac{G(z)}{Q(z)\left[Q(z)+z^2A(z)\right]}\\&=&1+z+\left(u_0+2\right)z^2+\left(b_3-u_0\right)z^3+\ldots
\end{eqnarray*}

De \cite{Dufresnoy1953} on obtient
\begin{equation*}
b_3-u_0\geq w_3
\end{equation*}
soit
\begin{equation*}
b_3\geq u_0^2+3u_0+2=\left(u_0+1\right)b_2
\end{equation*}

Avec la valeur $b_3=u_0^2+3u_0+2$ au développement de la fonction $B(z)/Q(z)$ à l’ordre 3 associons le polynôme
\begin{equation*}
D_4(z)=1-\left(u_0-\frac{3}{2u_0+3}\right)z+\frac{1}{2}(1+\frac{3}{2u_0+3})z^2+\left(u_0+1-\frac{3}{2u_0+3}\right)z^3-z^4
\end{equation*}

avec
\begin{equation*}
w_4=\left(u_0+1\right)b_3+u_0+\frac{1}{2}+\frac{3}{2u_0+3}\ ;D_4\left(u_0+1\right)=-u_0^2+{\frac{1}{2}u}_0+\frac{7}{4}
\end{equation*}
Comme précédemment :
\begin{equation*}
b_4-w_4\geq u_0
\end{equation*}

L’entier $b_4$ vérifie :
\begin{equation*}
b_4\geq\left(u_0+1\right)b_3+1
\end{equation*}

Avec la valeur
\begin{equation*}
b_4=\left(u_0+1\right)b_3+1
\end{equation*}

Le polynôme $D_5$ s’écrit
\begin{equation*}
D_5(z)=1-\left(u_0-\frac{1}{u_0+1}\right)z+\frac{1}{u_0+1}z^2+z^3+\left(u_0+1-\frac{1}{u_0+1}\right)z^4-z^5
\end{equation*}

Avec 
\begin{equation*}
w_5=\left(u_0+1\right)b_4-u_0+1+\frac{1}{u_0+1}\ ,\ D_5\left(u_0+1\right)=-\left(u_0^2-3\right)
\end{equation*}

L’entier $b_5$ vérifie
\begin{equation*}
b_5\geq w_5+u_0=\left(u_0+1\right)b_4+1+\frac{1}{u_0+1}
\end{equation*}
Soit
\begin{equation*}
b_5\geq\left(u_0+1\right)b_4+2
\end{equation*}

De l'équation (\ref{equ5}) on reprend la récurrence
\begin{equation*}
D_{n+2}(z)=\left(1+z\right)D_{n+1}(z)-z\left(\frac{b_{n+1}-w_{n+1}}{b_n-w_n}\right)D_n(z)\ ,\ n\geq1
\end{equation*}

Qu’on exploite pour $n=4$
\begin{equation*}
D_6\left(u_0+1\right)=-\left(u_0+2\right)\left(u_0^2-3\right)+\left(u_0+1\right)\left[\frac{b_5-w_5}{b_4-w_4}\right]\left[2u_0^2+\frac{3}{2}u_0-\frac{5}{4}\right]
\end{equation*}

Soit
\begin{equation*}
D_6(u_0+1)\geq u_0^2
\end{equation*}

Et par suite
\begin{equation*}
\theta>u_0+1
\end{equation*}
Inégalité qui est à rejeter. La valeur $b_2=u_0+2$ est à exclure. Il reste le développement
\begin{equation*}
\frac{B(z)}{Q(z)}=1+z+\left(u_0+1\right)z^2+\ldots
\end{equation*}

Et on conclut, pour $0\le d\le s-3$, par le développement
\begin{equation*}
f(z)=\frac{B(z)+zP(z)}{Q(z)+zA(z)}=1+z+z^2+\ldots
\end{equation*}

De ce développement, à l’ordre 2, de la fonction $f$, commun aux deux cas 
$d=s-2$ et  $d\le s-3$, on établit :
\begin{lemme}\label{lemme4} \hfill\\
$ u_0\geq2$, $d=s-2$ ou $d\le s-3$ et pour tout élément de l’ensemble $S^\prime(u_0)$ de fonction de rang infini $A(z)/Q(z)$, il existe un polynôme $U$, à coefficient entiers, tel que
\begin{itemize}
\item[1/]
\begin{equation}
\label{equ17}
\left\{\begin{array}{c}B(z)+zP(z)=\left[1-z^2-\varepsilon_2z^k\left(1+z+z^2\right)\right]U(z) \\
Q(z)+zA(z)=\left[1-z-z^2-\varepsilon_2z^k\left(1-z^2\right)\right]U(z)\end{array}\right. k\geq1
\end{equation}
\item[2/]  Le polynôme $U$ admet tous ses zéros sur $\left|z\right|=1$ et tous distincts vérifiant
\begin{equation*}
U(z)=\varepsilon_3 z^{d_1} U\left(\frac{1}{z}\right),\; \varepsilon_2\varepsilon_3=\varepsilon,\; k=s-d_1-1
\end{equation*}
\end{itemize}
\end{lemme}
A la fonction $f$ de rang fini et de développement
\begin{equation*}
f(z)=1+z+z^2+bz^3+\ldots
\end{equation*}

et à celle de rang infini et développement
\begin{equation*}
\frac{1-z^2}{1-z-z^2}=1+z+z^2+2z^3+\ldots
\end{equation*}

On associe la fonction
\begin{equation*}
g(z)=\frac{\left(1-z-z^2\right)f(z)-(1-z^2)}{\left(1-z^2\right)f(z)-(1+z-z^2)}
\end{equation*}
de développement
\begin{equation*}
g(z)=\left(b-2\right)\frac{z^3}{z^2}+\ldots
\end{equation*}

Les fonctions $f$ et $g$ sont de module 1 sur $\left|z\right|=1$.\newline
Du théorème de Rouché, la fonction $g$ admet au plus deux pôles dans $\left|z\right|<1$. Elle en admet exactement deux sur $\left|z\right|=1$.\newline
La fonction $g$ ne peut être identiquement nulle du fait que la fonction $f$ est de rang fini alors que celle 
\begin{equation*}
\frac{1-z^2}{1-z-z^2}
\end{equation*}
est de rang infini. La fonction $g$ se prolonge en une fonction holomorphe dans $\left|z\right|<1$, nulle à l’origine, à développement à coefficients entiers, de module égal à 1 sur $\left|z\right|=1$.\newline
Par le lemme de Schwarz
\begin{equation*}
g(z)=\varepsilon_2z^k\ ,\ \varepsilon_2=\pm1\ ,\ k\geq1
\end{equation*}

De cette expression de $g$ on déduit celle de $f$, à savoir 
\begin{equation*}
f(z)=\frac{B(z)+zP(z)}{Q(z)+zA(z)}=\frac{1-z^2-\varepsilon_2z^k(1+z-z^2)}{1-z-z^2-\varepsilon_2z^k(1-z^2)}
\end{equation*}
La formation du système (\ref{equ17}) en découle avec la présence d’un polynôme $U(z)$.
Le polynôme
\begin{equation*}
Q(z)+zA(z)
\end{equation*}

admet exactement un zéro dans  $\left|z\right|<1$. De la valeur $U(0)=1$ et du lemme 3 de \cite{Dufresnoy1953} le polynôme $U$ admet tous ses zéros sur $\left|z\right|=1$ et tous distincts. Pour accompagner le système (\ref{equ17}) on reprend de (\ref{equ14}) le polynôme $V$ de degré $s-2$ et de (\ref{equ7}) le polynôme $V$ de degré  $\le s-3$. 

On établit :
\begin{lemme}\label{lemme5}\hfill
\begin{itemize}
\item[1/] Pour $u_0\geq2$ , $ s=h$ , $ d=s-2$\newline
\qquad i) Les polynômes $U$ et $V$ vérifient
\begin{equation*}
U(z)=1+z\ ;V(z)=1-\varepsilon z^{s-2}\ ,\ \varepsilon=\pm1\ ,\ s\geq3
\end{equation*}
 ii)	Seule la suite ($\gamma_s$) représentée par la famille de fonctions de rang infini où
\begin{equation*}
\left.\begin{array}{l}A(z) =u_0+\left(u_0-2\right)z-z^2-\varepsilon z^{s-2}\left[u_0+\left(u_0-1\right)z-z^2\right] \\
Q(z)=1-u_0z-u_0z^2-\varepsilon z^{s-2}\left[1-\left(u_0-1\right)z-u_0z^2\right]\end{array}\right. s\geq3 
\end{equation*}
\begin{equation*}
G(z) =\left[\left(1-z\right)\left(1-\varepsilon z^{s-2}\right)\right]^2
\end{equation*}
appartient à l’ensemble $S^\prime(u_0)$
\item[2/] Pour $u_0\geq2$, $ d=1=s-3$\newline
\qquad i) Les polynômes $U$ et $V$ vérifient
\begin{equation*}
U(z)=1-z^2\ ,\ V(z)=1+z
\end{equation*}
 ii)	Seul l’élément $\gamma_4$, $\varepsilon=$1, représenté par la fonction de rang infini où
\begin{align*}
A(z)&=u_0+\left(u_0-1\right)z-u_0z^2-\left(u_0-1\right)z^3+z^4\\
Q(z)&=1-u_0z-\left(u_0+1\right)z^2+\left(u_0-1\right)z^3+u_0z^4\\
G(z)&=\left[(1+z)(1-z^2\right]^2
\end{align*}
appartient à l’ensemble $S^\prime(u_0)$.
\item[3/] Pour $u_0\geq2$, $s=h$, $0\le d\le s-3$\newline
\qquad i) Les polynômes $U$ et $V$ vérifient
\begin{equation*}
U(z)=V(z)=1
\end{equation*}
ii)	Seule la suite ($\alpha_s$), représentée par la famille de fonctions de rang infini où
\begin{equation*}
\left.\begin{array}{l}A(z)=u_0-z+u_0z^{s-1}-z^s \\
Q(z)=1-\left(u_0+1\right)z+z^{s-1}-u_0z^s\end{array}\right. s\geq3 
\end{equation*}
\begin{equation*}
G(z)=A(z)-z^{s-2}Q(z)\end{equation*}
appartient à l’ensemble $S^\prime(u_0)$.
\end{itemize}
\end{lemme}
La suite ($\alpha_s$) peut être prolongée à $s=2$ par 
\begin{align*}
A(z)&=u_0+\left(u_0-1\right)z-z^2\\
Q(z)&=1-u_0z-u_0z^2\\
G(z)&=A(z)-Q(z)
\end{align*}

\paragraph{- Etude du cas $d=s-2$}\hfill\\
Des relations (\ref{equ14}) et (\ref{equ17}) on établit l’identité
\begin{equation*}
(1+z)V(z)=\left(1-\varepsilon z^{s-d_1-1}\right)U(z)
\end{equation*}
de laquelle on tire
\begin{equation*}
U(z)=1+z\ ,\ V(z)=1-\varepsilon z^{s-2}\ ,\ s\geq3
\end{equation*}
et de (\ref{equ16})
\begin{equation*}
A(z)-Q(z)=\left(u_0-1\right)\left(1+z\right)^2(1-\varepsilon z^{s-2})
\end{equation*}
et de (\ref{equ17})
\begin{equation*}
Q(z)+zA(z)=(1+z)\left[1-z-z^2-\varepsilon z^{s-2}\left(1-z^2\right)\right]
\end{equation*}
et on conclut :
\begin{equation}
\label{equ18}
\left.\begin{array}{l}A(z)=u_0+\left(u_0-2\right)z-z^2-\varepsilon z^{s-2}(u_0+\left(u_0-1\right)z-z^2) \\Q(z)=1-u_0z-u_0z^2-\varepsilon z^{s-2}(1-\left(u_0-1\right)z-u_0z^2)\end{array}\right. s\geq3
\end{equation}

\paragraph{- Etude du cas $d=1=s-3$}\hfill\\
Des relations (\ref{equ9}) et (\ref{equ17}) on établit l’identité
\begin{equation}
\label{equ19}
\left(1-\varepsilon\varepsilon_1z^{s-d-2}\right)V(z)=\left(1-\varepsilon\varepsilon_3z^{s-d_1-2}\right)U(z) 
\end{equation}

De laquelle on tire
\begin{equation*}
V(z)=1+z\ ,\ U(z)=1-z^2 \text{ avec } d_1=2\ ,\ s=4\ ,\ \varepsilon_3=-1
\end{equation*}

De l’expression de $U(z)=1-z^2=1-z^{4-2}$,  $ \varepsilon=1$, on retrouve le polynôme $Q$ de degré 4 de (\ref{equ19})
\begin{equation*}
Q(z)=1-u_0z-\left(u_0+1\right)z^2-\left(u_0-1\right)z^3+z^4
\end{equation*}

Et de (\ref{equ17}) on complète par 
\begin{equation*}
A(z)=u_0+\left(u_0-1\right)z-u_0z^4-\left(u_0-1\right)z^3+z^4
\end{equation*}

\paragraph{- Etude du cas $0\le d\le s-3$}\hfill\\
La relation (\ref{equ19}) est aussi vérifiée par
\begin{equation*}
U(z)=V(z)=1
\end{equation*}

Et on retrouve le résultat de Dufresnoy-Pisot dans \cite{Dufresnoy1953}.
Parmi les éléments de degré $s=h=3$ de l’ensemble $S^\prime(u_0)$ seul l’élément $\alpha_3$ est représenté par la fonction de rang infini où
\begin{eqnarray*}
A(z)=u_0-z+u_0z^2-z^3\ ;Q(z)=1-\left(u_0+1\right)z+z^2-u_0z^3\\
G(z)=A(z)-zQ(z);\frac{Inf}{\left|z\right|=1}.\frac{G(z)}{z^2}=1-\frac{1}{u_0}
\end{eqnarray*}

De plus pour $u_0=1$ on retrouve l’élément ($\alpha_3$) de la suite ($\alpha_s$) représentée par la famille de fonctions de rang infini, Chapitre III de \cite{Amara1966}
\begin{equation*}
A(z)=1-z+z^{s-1}-z^s\ ;Q(z)=1-2z+z^{s-1}-z^s
\end{equation*}
Par un procédé heuristique seule la suite ($\alpha_s$) représentée par la famille de fonctions de rang infini où
\begin{align*}
A(z)&=u_0-z+u_0z^{s-1}-z^s\ ;Q(z)=1-\left(u_0+1\right)z+z^{s-1}-u_0z^s \; s\geq3 \\
G(z)&=A(z)-z^{s-2}Q(z)
\end{align*}
appartient à l’ensemble $S^\prime(u_0)$.

En regroupant les résultats des lemmes \ref{lemme1} et \ref{lemme5}, on termine la première partie de notre étude par :
\begin{theoreme} \label{th1}\hfill\\
1/ Pour $u_0\geq2$, $s\geq h+1$ seule la suite ($\beta_s$) représentée par la famille de fonction de rang infini où
\begin{equation*}
\left.\begin{array}{l}\left(1-z\right)A(z)=u_0-z^{s-1}-\left(u_0-1\right)z^s \\
\left(1-z\right)Q(z)=1-\left(u_0+1\right)z+u_0z^{s+1}
\end{array}\right.
s\geq2
\end{equation*}
appartient à l’ensemble $S^\prime(u_0)$. De plus cette famille vérifie
\begin{align}
\label{equ20}
Q(z)&+zA(z)=P(z)-B(z)=1-z^s\nonumber \\
Q(z)&+zP(z)=1-z^{s+1}\\
G(z)&=u_0\left(1-z^s\right)^2\nonumber
\end{align}
2/ Pour $u_0\geq2$, $d=s-2\geq1$, seule la suite ($\gamma_s$) représentée par la famille de fonctions de rang infini où
\begin{equation*}
\left.\begin{array}{l}
{A(z)=u}_0+\left(u_0-2\right)z-z^2-\varepsilon z^{s-2}(u_0+\left(u_0-1\right)z-z^2)\\
Q(z)=1-u_0z-u_0z^2-\varepsilon z^{s-2}(1-\left(u_0-1\right)z-{u_0z}^2)
\end{array}\right.
\varepsilon=\pm1,\ s\geq3
\end{equation*}
appartient à l’ensemble $S^\prime(u_0)$. De plus, cette famille vérifie 
\begin{align}
\label{equ21}
Q(z)&+zA(z)=\left(1+z\right)\left[1-z-z^2-\varepsilon z^{s-2}\left(1-z^2\right)\right] \nonumber \\
Q(z)&+zP(z)=1-z^2-z^3-\varepsilon z^{s-2}\left(1+z-z^3\right)\\
G(z)&=\left(u_0-1\right)\left[\left(1+z\right)\left(1-\varepsilon z^{s-2}\right)\right]^2\nonumber
\end{align}
3/ Pour $u_0\geq2$, $d=1=s-3$, seul l’élément $\gamma_4$,  $\varepsilon=1$, représenté par la fonction de rang infini où
\begin{align*}
A(z)&=u_0+\left(u_0-1\right)z-u_0z^2-\left(u_0-1\right)z^3+z^4\\
Q(z)&=1-u_0z-\left(u_0+1\right)z^2-\left(u_0-1\right)z^3+z^4
\end{align*}
appartient à l’ensemble $S^\prime(u_0)$. De plus cet élément vérifie 
\begin{align}
\label{equ22}
Q(z)&+zA(z)=(1+z) \left[1-z-z^2-z^2 (1-z^2)\right] \nonumber \\
Q(z)&+zP(z)= (1+z) \left[1-z-z^2-z^3+ z^4 \right] \\
G(z)&=u_0\left[\left(1-z\right)\left(1-z^2\right)\right]^2\nonumber
\end{align}
4/ Pour $u_0\geq2$, $d=0\le s-2$, seule la suite ($\alpha_s$) représentée par la famille de fonctions de rang infini où
\begin{align*}
A(z)&=u_0-z+u_0z^{s-1}-z^s \\
Q(z)&=1-\left(u_0+1\right)z+z^{s-1}-u_0z^s
\end{align*}
appartient à l’ensemble $S^\prime(u_0)$. De plus cette famille vérifie
\begin{align}
\label{equ23}
Q(z)&+zA(z)=1-z-z^2+z^{s-1}(1-z^2) \nonumber \\
Q(z)&+zP(z)=1-z-z^2+z^{s-1}(1+z-z^2)\\
G(z)&=A(z)-z^{s-2}Q(z)\nonumber
\end{align}
\end{theoreme}
\hfill\\
Remarques sur le théorème \ref{th1} : \newline
1/ Le cas $u_0=1$, représenté par l’ensemble $S^\prime(1)$, avait fait l’objet du chapitre III de la thèse de l’auteur \cite{Amara1966}.\newline
2/ Pour $u_0=1$, la suite ($\gamma_s$), $\varepsilon=\pm1$, cesse d’appartenir à l’ensemble $S^\prime$.\newline
3/ Un résultat intéressant réside dans l’indépendance des polynômes
\begin{equation*}
Q(z)+zA(z)\ ,\ Q(z)+zP(z)
\end{equation*}
du coefficient $u_0\geq2$.\newline
4/ Signalons l’apparition des polynômes des équations (\ref{equ21}) et (\ref{equ23})
\[ i)\; Q(z)+zP(z)=1-z^2-z^3-\varepsilon z^{s-2}(1+z-z^3)\]
\[ ii)\; Q(z)+zP(z)=1-z-z^2+z^{s-1}\left(1+z-z^2\right)\]

La famille $ii)$ représente, pour la valeur de $s$ convenablement choisie, une suite de nombre de Salem convergente par valeurs inférieures vers l’élément $\theta^\prime$, plus petit élément de $S^\prime$, zéro du polynôme $1+z-z^2$.\newline
Les deux familles de $i)$ représentent, pour des valeurs de $s$ convenablement choisies, deux suites de nombres de Salem, la première, pour $\varepsilon=1$, convergente par valeurs supérieures, la deuxième pour $\varepsilon=-1$, convergente par valeurs inférieures vers l’élément $\theta_0$, plus petite élément de l’ensemble $S$, zéro du polynôme
\begin{equation*}
1+z-z^3
\end{equation*}
Le plus petit élément $\tau_0$ de la suite convergente par valeurs inférieures vers l’élément $\theta_0$ est le zéro du polynôme issu de (\ref{equ21}) avec $s=10$, $\varepsilon=-1$,
\begin{equation*}
1+z-z^3\left(1+z+z^2+z^3+z^4\right)+z^9+z^{10}
\end{equation*}
qui n’est autre que le plus petit élément connu des nombres de Salem. 

\section{Sur les points d’accumulation des nombres de Salem}
\subsection{Définitions} On entend par nombre de Salem un entier algébrique $\tau$ supérieur à 1 dont les conjugués, tous distincts et à l’exception des nombres $\tau$ et $\frac{1}{\tau}$, sont situés, au moins au nombre de deux, sur $\left|z\right|=1$. Le polynôme à coefficients entiers, dont le nombre $\tau$ est zéro, est désigné par la lettre $R$ s’écrivant
\begin{equation*}
R(z)=\left[1-\left(\tau+\frac{1}{\tau}\right)z+z^2\right]\prod_{1}^{s-1}\left[1-\left(\alpha_{i\ }+\bar{\alpha_i}\right)z+z^2\right]
\end{equation*}
avec degré de $R=2s$, $s\geq2$,  $\alpha_i\neq\bar{\alpha_i}$ et $\left|\alpha_i\right|=1$.\newline
Le polynôme $R$ est à la fois réciproque et irréductible.\newline
L’ensemble des nombres de Salem est désigné par la lettre $T$. Salem a établi, dans \cite{Salem1944}, un lien étroit entre les deux ensembles $S$ et $T$ à savoir :
\begin{lemme}\label{lemme5bis}
Tout élément de l’ensemble S est un point d’accumulation de l’ensemble T.
\end{lemme}
\begin{itemize}
\item[$\bullet$]	A cet effet à tout élément $\theta$ de l’ensemble $S$, représenté par sa fonction de rang fini $P(z)/Q(z)$,  $P\neq Q$ et à tout entier $m\geq1$ et  $\eta=\pm1$, il associe un nombre $\tau_m$ de l’ensemble $T$ représenté par son polynôme $R_m$ vérifiant la relation
\begin{equation}
\label{equ24}
U_m(z)R_m(z)=Q(z)+\eta z^mP(z)
\end{equation}
avec degré de $P=s$, $P(0)=u_0\geq1$ et fixé.
\end{itemize}
Le polynôme $U_m$, à zéros tous situés sur $\left|z\right|=1$, est choisi rendant irréductible le polynôme $R_m$. Dufresnoy et Pisot ont dans \cite{Dufresnoy1953} établi que le polynôme du second membre de (\ref{equ24}) admet au plus un zéro dans $\left|z\right|<1$ et au moins $(m+s-2)$ zéros distincts sur $\left|z\right|=1$. La suite ($\tau_m$), zéros des polynômes $R_m$ à partir d’un certain rang $m$ appartient à l’ensemble $T$ et vérifie
\begin{equation*}
\lim_{m\rightarrow+\infty}{\tau_m=\theta}
\end{equation*}
Dans l’étude qui va suivre on propose une réciproque au procédé de Salem \cite{Salem1944}. A cet effet on utilise le procédé de Boyd \cite{Boyd1977}, réciproque de celui de Salem à savoir :\newline
A tout élément $\tau$ de l'ensemble $T$, zéro d'un polynôme $R$, on associe un élément $\theta$ de l'ensemble $S$, représenté pas sa fonction de rang fini $P(z)/Q(z)$, $P\neq Q$, $P(0)=u_0\geq1$ et l'entier $m\geq1$, $\eta=\pm1$ vérifiant la relation
\begin{equation}
\label{equ24'}
U(z).R(z)=Q(z)+\eta z^m P(z)
\end{equation}
\begin{lemme}
\label{lemme6}
 Des procédés de Salem et Boyd, on établit : Pour tout élément de l’ensemble $T$, associé à l’élément de l’ensemble $S$ par la relation (\ref{equ24'}).\newline
1/ Pour $m\geq2$ et $\eta=\pm1$ ou $m=1$ et $\eta=-1$, $u_0\geq1$, fixés.\newline
L’élément $\tau$ est borné si et seulement si l’élément $\theta$ est borné.\newline
2/ Pour $m=1$ et $\eta=+1$, $u_0\geq1$ .\newline
L’élément $\tau$ est borné si et seulement si l’expression $\left(\theta-u_0\right)$, strictement positive, est bornée.\newline
\end{lemme}

Considérons la représentation de la fonction
\begin{equation*}
\frac{P(z)}{Q(z)}=\frac{\varphi(z)}{\Psi(z)}
\end{equation*}
en un quotient de deux produits de Blaschke à savoir
\begin{equation*}
\Psi(z)=\frac{1-\theta z}{\theta-z}\ ;\ \varphi(z)=\prod_{1}^{s-1}{\left(\frac{1-\bar{\theta_i}z}{\theta_i-z}\right)\frac{\theta_i}{\left|\theta_i\right|}\ ;\ \left|\theta_i\right|>1\ ,\ \varphi(0)=\frac{u_0}{\theta}}
\end{equation*}
De la relation (\ref{equ24'}) formons la fonction
\begin{equation*}
g(z)=\left(1-\theta z\right)\frac{U_m(z)R_m(z)}{Q(z)}=1-\eta z^{m+1}\varphi(z)-\theta z\left[1-\eta z^{m-1}\varphi(z)\right]
\end{equation*}
La fonction $g$ est holomorphe et à zéro unique $\frac{1}{\tau}$ dans $\left|z\right|<1$.\newline

1/ Pour $m\geq2$ et $\eta=\pm1$ ou $m=1$ et $\eta=-1$ fixés\newline
La famille de fonctions $g$ constitue, comme l’est la famille de fonctions $\varphi$, une famille de fonctions holomorphes et bornées dans $\left|z\right|<1$ si et seulement si l’élément $\theta$ est borné.\newline
D’après Montel, de la famille de fonctions $g$, on peut extraire une suite de fonctions uniformément convergente sur tout compact de $\left|z\right|<1$. De la valeur $g(0)=1$, il existe un réel $\delta>0$ tel que sur le compact $\left|z\right|\le\delta$, la suite de fonctions $g$ ne s’annule pas et par conséquent l’élément $\tau$ est borné si et seulement si l’élément $\theta$ est borné.\newline

2/ Pour $m=1$ et $\eta=1$\newline
La fonction $g$ s’écrit
\begin{equation*}
g(z)=1-z^2\varphi(z)-\theta z\left[1-\varphi(z)\right]
\end{equation*}
Pour $u_0$ tendant vers $+\infty$, de l’encadrement $u_0<\theta<u_0+1$, les valeurs $\frac{u_0}{\theta}$ et $\varphi\left(\frac{u_0}{\theta}\right)$ tendent simultanément vers 1 et on se trouve devant l’expression non exploitable
\begin{equation*}
u_0\left[1-\varphi\left(\frac{u_0}{\theta}\right)\right]
\end{equation*}
Pour contourner cette difficulté on utilise la méthode Schur en formant la fonction
\begin{equation*}
z\varphi_1(z)=\frac{\varphi(z)-\frac{u_0}{\theta}}{1-\frac{u_0}{\theta}\varphi(z)}
\end{equation*}
La fonction $\varphi_1$ est aussi un produit de Blaschke permettant d’écrire
\begin{equation*}
\varphi(z)=\frac{\frac{u_0}{\theta}+z\varphi_1(z)}{1+\frac{u_0}{\theta}z\varphi_1(z)}
\end{equation*}
entrainant pour la fonction $g$ l’expression
\begin{equation*}
g(z)=1-z^2\varphi(z)-\left(\theta-u_0\right)z\frac{1-z\varphi_1(z)}{1+\frac{u_0}{\theta}z\varphi_1(z)}
\end{equation*}
Par le même raisonnement que précédemment l’élément $\tau$ est borné si et seulement si l’expression $\left(\theta-u_0\right)$ est bornée.\newline
Du lemme \ref{lemme6}, paragraphe 1,\newline
- Pour $m\geqslant2$ et $\eta=\pm1$ ou $m=1$ et $\eta=-1$ et pout tout $u_0$ fixé, $u_0\geqslant1$, des suites ($\tau_\nu$) et ($\theta_\nu$) simultanément bornées avec $\tau_\nu < u_0 +1$, on peut extraire des sous-suites convergentes et tout en gardant les mêmes notations,
\begin{equation*}
\lim_{\nu\rightarrow\infty}{\theta_\nu=\theta} \text{ et, pour tout } \left|z\right|<1 ,\;\lim\limits_{\nu\rightarrow\infty} \frac{P_\nu(z)}{Q_\nu(z)}=\frac{A(z)}{Q(z)}
\end{equation*}
avec l'élément $\delta$ appartenant à l'ensemble $S^{\prime}$ et représenté par la fonction de rang infini $A(z)/Q(z)$. \newline
\begin{equation*}
 \lim_{\nu\rightarrow\infty}{\tau_v=\tau^\prime} \text{ et, pour tout} \left|z\right|<1 ,\;
\lim_{\nu\rightarrow\infty}\frac{U_\nu(z) R_\nu(z)}{Q_\nu(z)}=\frac{Q(z)+\eta z^m A(z)}{Q(z)}
\end{equation*}
avec l'élément $1/\tau^{\prime}$ zéro du polynôme $ Q(z)+\eta z^mA(z) $.\newline
De la valeur $P_\nu(0)= u_0$, on dispose de celle $A(0)=u_0$ et des fonctions
\begin{equation*}
\frac{A(z)}{Q(z)}\frac{1-\theta z}{\theta-z} ,\;\frac{P(z)}{Q(z)}\frac{1-\theta z}{\theta-z}
\end{equation*}
holomorphes et bornées dans $|z| <1$. On en déduit,
\[ u_0 < \theta< u_0 +1,\; P(0) = v_0 < \theta < u_0 + 1 \]
soit \[ P(0) = v_0 \leqslant u_0  \]
Etudions l'élément $\tau^\prime$ qui présente deux possibilités : ou bien il appartient à l'ensemble $S$ ou bien il est égal à 1.\newline
Par le théorème de Rouché, le polynôme 
\[ Q(z) + \eta z^m A(z)\]
admet au plus un zéro dans $\left|z\right|<1$. Si ce zéro appartenait à $\left|z\right|<1$, l'élément $\tau^\prime$ serait bien dans l'ensemble $S$.\newline
Si ce polynôme admettait tous ces zéros dans $\left|z\right|\geqslant1$, on aurait $1/\tau^{\prime}\; \geqslant 1$, soit $\tau^{\prime} \leqslant 1$ et, par suite, $\tau^{\prime} = 1$.\newline
Examinons le cas $P(0)=v_0 = u_0$, soit $P(0)=A(0) = u_0$. L'élément $\theta$ vérifiant l'encadrement 
\[ u_0 < \theta < u_0 + 1\]
appartient à l'ensemble $S^\prime (u_0)$. Des compléments sur les contenus du théorème \ref{th1} (paragraphe 4), la famille de polynômes 
\[  Q(z) + z P(z) =1-z^2-z^3-\varepsilon z^{s-2}(1+z-z^3) \]
formée de la représentation de l'élément $\gamma_4$ (\ref{equ22}) et de celle de la suite ($\gamma_s$)  (\ref{equ21}) de rang infini si et seulement si $u_0 \geqslant 2$.\newline
- Du lemme \ref{lemme4}, paragraphe 2, pour $m = 1$,  $\eta = + 1$ et pour tout $u_0$ fixé et $u_0 \geqslant 1$ \newline
A toute suite ($\tau_{\nu}$) bornée, on associe la suite ($\theta_{\nu}$) vérifiant
\[ 0 < \theta_{\nu} - u_0 < 1\]
Des suites bornées ($\tau_{\nu}$) et ($\theta_{\nu}$), on peut extraire des sous-suites convergentes et, tout en gardant les mêmes notations,
\begin{align}
\label{equ25}
\text{1/ }& \lim_{\nu\rightarrow\infty}{\theta_\nu=\theta} \text{et pour } \left|z\right|<1 ,\;\lim\limits_{\nu\rightarrow\infty} \frac{P_\nu(z)}{Q_\nu(z)}=\frac{A(z)}{Q(z)}\\
\text{2/ }& \lim_{\nu\rightarrow\infty}{\tau_v=\tau^\prime} \text{et pour } \left|z\right|<1 ,\;
\lim_{\nu\rightarrow\infty}\frac{U_\nu{(z) R}_\nu(z)}{Q_\nu(z)}=\frac{Q(z)+ z A(z)}{Q(z)}\nonumber
\end{align}
Du lemme \ref{lemme4}, paragraphe 1, \newline
- l'élément $\tau^\prime$ présente deux possibilités :
\begin{equation}
\label{equ26} \text{ou bien il appartient à l'ensemble $S$ ou bien il est égal à 1.}
\end{equation}
- de la valeur $P(0) = v_0 =u_0 = A(0)$ associée à l'encadrement
\[u_0 < \theta < u_0 + 1\]
l'élément $\theta$ appartient à $S^\prime (u_0)$. Pour $u_0 \geqslant 2$, on retrouve la famille de polynômes de (\ref{equ17}).\newline
Pour la suite de notre étude, il nous reste à reprendre de \cite{Amara1966} la composition de l'ensemble $S^\prime (1)$ et à déterminer les propriétés de l'ensemble 
\[ S^\prime (u_0 , v_0) = S^{\prime} \cap ] u_0, u_0 + 1 [ \;, \;A(0) = u_0 \;, \; 1\leqslant P(0) = v_0 \leqslant u_0 -1 \]

\begin{theoreme} \label{th2}\hfill\\
 1/ La suite ($\beta_s$) représenté par les deux familles de fonctions de rang infini
\begin{align*}
\left(1-z\right)&A_1(z)=1-z^{s-1}\ ;\left(1-z\right)A_2(z)=1-z-z^s+z^{s+1}\\
\left(1-z\right)&Q(z)=1-2z+z^{s-1}
\end{align*}
appartient à l’ensemble $S^\prime(1)$\newline
2/ Pour $s=h=4$, l’élément $\gamma_4$ représenté par la fonction de rang infini
\begin{equation*}
A(z)=1-z^2+z^4\ ;Q(z)=1-z-2z^2+z^4
\end{equation*}
appartient à l’ensemble $S^\prime(1)$\newline
3/ La suite ($\alpha_s$) représentée par les deux familles de fonctions de rang infini
\begin{align*}
A_1(z)&=1-z-z^{s-1}\ ;A_2(z)=1-z+z^{s-1}-z^s\\
Q(z)&=1-2z+z^{s-1}-z^s
\end{align*}
appartient à l’ensemble $S^\prime(1)$.
\end{theoreme}

\subsection{Sur la composition de l’ensemble $S^\prime\left(u_0,\ v_0 \right)$}
\subsubsection{Plus petit élément de $S^\prime(u_0, v_0)$}
\paragraph{- Etude du cas $h\geq s+1$}\hfill\\
De $G(z)=A(z)B(z)-\varepsilon\varepsilon^\prime z^{h-s}P(z)Q(z)$, on obtient $G(0)=u_0b_0$.
De l'équation (\ref{equ1}) de la première partie on reprend 
\[\theta^2\geq v_0a_1=u_0f_1(0)=u_0^2+u_0b_2+b_0^2+u_0b_0\geq\left(u_0+1\right)^2\]
Soit  $\theta\geq u_0+1$, une inégalité à rejeter. \newline
\paragraph{- Etude du cas $s\geq h+1$}\hfill\\
De $G(z)=P(z)Q(z)-\varepsilon\varepsilon^\prime z^{s-h}A(z)B(z)$, on obtient $G(0)=v_0$.
De l'équation (\ref{equ2}), on tire
\[u_0a_1=u_0^2+2v_0\]
De (\ref{equ4}) on déduit
\[\theta^4-\left(u_0^2+2u_0\right)\theta^2+v_0^2\geq0\]
	ou
\[\theta^2-\theta u_0-v_0\ \geq0\]

Le plus petit élément de $S^\prime(u_0, v_0)$, pour $s\geq h+1$, est représenté par la fonction de rang infini
\[A(z)=u_0+\left(v_0-1\right)z\ ;Q(z)=1-u_0z-v_0z^2\]
	ou, autrement
\begin{align}
\label{equ27}
\left(1-z\right)A(z)&=u_0-\left(u_0-v_0+1\right)z+\left(v_0-1\right)z^2\\
\left(1-z\right)Q(z)&=1-\left(u_0+1\right)z+\left(u_0-v_0\right)z^2+v_0z^3\nonumber
\end{align}
\newline
\paragraph{- Etude du cas $s=h$}\hfill\\
De $G(z)=A(z)B(z)\ -\ P(z)Q(z)$, $G(0)=u_0b_0-v_0$ et par suite $k=0$ 
on obtient de (\ref{equ1})
\[v_0f_1(0)=u_0^2+b_0^2+2G(0)\geq u_0^2+2\left(u_0-v_0\right)+1\ ,\ b_0=1\]

	et on retrouve (\ref{equ7})

\[\theta^4-\left[u_0^2+2\left(u_0-v_0\right)+1\right]\theta^2+v_0^2\geq0\]
	ou
\[\theta^2-\left(u_0+1\right)\theta+v_0\geq0\]
Le plus petit élément de $S^\prime(u_0, v_0)$, pour $s=h$, est représenté par la fonction de rang infini
\[A(z)=u_0-\left(v_0+1\right)z+z^2\ ;Q(z)=1-\left(u_0+1\right)z+v_0z^2\]
En s’inspirant de la représentation des plus petits éléments de $S^\prime(u_0, v_0)$, on établit la composition de l’ensemble $S^\prime(u_0, v_0)$.

\begin{theoreme} \label{th3}\hfill\\
1/ Pour $u_0\geq2\ ,\ 1\le v_0\le u_0-1$, $s=h+1$, seule la suite ($\beta_s$) représentée par la famille de fonctions de rang infini où
\begin{align*}
\left(1-z\right)A(z)&=u_0-\left(u_0-v_0+1\right)z^{s-1}-\left(v_0-1\right)z^s\\
\left(1-z\right)Q(z)&=1-\left(u_0+1\right) z +\left(u_0-v_0\right) z^s + v_0\; z^{s+1}\; , s\geq2
\end{align*}
appartient à l’ensemble $S^\prime(u_0, v_0)$. De plus cette famille vérifie :
\begin{align*}
Q(z)+zA(z)&=P(z)-B(z)=1-z^s\\
G(z)&=v_0\left(1-z^s\right)^2
\end{align*}
2/ Pour $u_0\geq2$, $ 1\le v_0\le u_0-1$, $s=h$\newline
Chaque élément de $S^\prime(u_0, v_0)$ de fonction de rang infini $A(z)/Q(z)$ représente l’une des deux possibilités de (\ref{equ26}) :
\begin{itemize}
\item[i)] Le polynôme
\[Q(z)+\eta \; z^mA(z)\]
admet tous ces zéros sur $\left|z\right|=1$, y compris l’élément $\tau^\prime=1$.

Dans ce cas seul le plus petit élément de $S^\prime(u_0, v_0)$ représenté par la fonction de rang infini où
\[A(z)=u_0-\left(v_0+1\right)z+z^2\ ;Q(z)=1-\left(u_0+1\right)z+v_0\ z^2\]
et vérifiant
\begin{align*}
Q(z)+zA(z)&=B(z)+zP(z)=\left(1-z\right)^2(1+z)\\
G(z)&=\left(u_0-v_0\right)\left(1-z^2\right)^2
\end{align*}
y répond avec $m=1$, $\eta=+1$.\newline
\item[ii)]  Pour $m\geq2$, $\eta=\pm1$, $m=1$, $\eta=-1$, \newline
Le polynôme 
\[Q(z)+\eta z^m A(z)\]
Admet exactement un zéro $1/\tau^\prime$ dans $\left|z\right|<1$ avec le point $\tau^\prime$ appartenant à l’ensemble $S$.
\end{itemize}
\end{theoreme}
\hfill\\
1/ Du transfert de la valeur $P(0)=u_0$ à celle de $P(0)=v_0$ par un procédé heuristique sur les résultats du théorème \ref{th1}, pour $s=h+1$, on obtient la première partie du théorème \ref{th3}\newline
2/ Examinons seulement le cas des polynômes
\[Q(z)+\eta z^mA(z)\]
admettant tous leurs zéros sur $\left|z\right|=1$. De la valeur $Q(0)=1$, ces polynômes sont réciproques, et par suite on obtient
\[B(z)+\eta\varepsilon\varepsilon^\prime z^mP(z)=Q(z)+\eta z^mA(z)\]
	Ou mieux
\[B(z)-Q(z)=\eta z^m\left[A(z)-\varepsilon\varepsilon^\prime P(z)\right]\]
De la valeur $G(0)=u_0-v_0$ on tire $\varepsilon\varepsilon^\prime=1$. Du théorème de Rouché on a $m=1$ ; du développement de la fonction
\[\frac{B(z)}{Q(z)}=1+b_1z+\ldots\]
avec $b_1>0$ on a nécessairement
\[b_1=u_0-v_0  \text{ et }  \eta=1\]	
On termine par la relation
\[B(z)-Q(z)=z\left[A(z)-P(z)\right]\]
De laquelle on forme le système
\begin{align*}
B(z)-Q(z)&=zU(z)\\
A(z)-P(z)&=U(z)
\end{align*}

Des résultats précédents le polynôme $U$ admet tous ses zéros sur $\left|z\right|=1$ et tous distincts. Il vérifie en plus
\[U=\varepsilon z^{s-1}U\left(\frac{1}{z}\right)\ ,\ U(0)=u_0-v_0\]
Du système précédent on déduit
\[G(z)=\left[Q(z)+zA(z)\right]U(z)\]

Du degré de $G$ égal à $2s$, de celui de $U$ égal à $s-1$, de l’unicité des zéros de $U$, de $U(0)=u_0-v_0$, il existe un entier $k$ tel que
\[G(z)=\left(u_0-v_0\right)\left[1-kz+z^2\right]U_1^2(z)\]
	Avec $U(z)=\left(u_0-v_0\right)U_1(z)$, $U_1(0)=1$
De l’écriture de $G$ sur $\left|z\right|=1$
\[G(z)=-\varepsilon z^s\left[\left|Qz)\right|^2-\left|A(z)\right|^2\right]=-\varepsilon\left[k-\left(z+\frac{1}{z}\right)\right]z^s\left|U_1(z)\right|^2\]
On déduit
\[k\geq2\]
Posons
\[\frac{U_1(z)}{Q(z)}=1+a_1z+\ldots\]
De l’expression du polynôme $G$
\[1+z\frac{A(z)}{Q(z)}=(1-kz+z^2)\frac{U_1(z)}{Q(z)}\]
On obtient
\[a_1=u_0+k\]
Au développement de 
\[\frac{B(z)}{Q(z)}=1+\left(u_0-v_0\right)z\frac{U_1(z)}{Q(z)}=1+\left(u_0-v_0\right)z+\left(u_0-v_0\right)\left(u_0+k\right)z^2+\ldots\]
On associe le polynôme
\[D_3(z)=1-\left(v_0+k-1\right)z+\left(u_0+k-1\right)z^2-z^3\]

Pour $k\geq3$ le polynôme $D_3$ vérifie
\[D_3\left(u_0+1\right)>0\]
L’élément $\theta$ de $S^\prime(u_0,v_0)$ est alors strictement supérieur à $\left(u_0+1\right)$. Ce résultat est à rejeter, soit
$k\le2$ et, par suite, $k=2$.
Il en résulte \[Q(z)+zA(z)=B(z)+zP(z)=\left(1-z\right)^2U_1(z)\]
On retrouve le plus petit élément de $S^\prime(u_0,v_0)$ pour $s=h$, avec $U_1=1+z$.
En s’inspirant de la représentation du plus petit élément de $S^\prime (u_0,v_0)$, pour $s=h$, et en le désignant par la lettre $\alpha_2$, on forme la suite ($\alpha_s$) appartenant à l’ensemble $S^\prime$ et représenté par la famille de fonctions de rang infini où,
\begin{align*}
A(z)=u_0+z-\left(v_0+2\right)z^{s-1}+z^s \\Q(z)=1-\left(u_0+2\right)z+z^{s-1}+v_0z^s\\
G(z)=\left(u_0-v_0\right)\left[(1-z)(1+z^{s-1}\right]^2
\end{align*}
Toutefois, pour $s\geq3$, on obtient
\[ \alpha_s>u_0+1 \]
De plus cette famille vérifie
\begin{eqnarray}
\label{equ28}
\\
Q(z)+z A(z) = B(z)+zP(z) = B(z)+zP(z)=\left(1-z\right)^2\left(1+z^{s-1}\right)  \nonumber \\ 
Q(z)-zP(z) =1-\left(u_0+v_0+2\right)z-z^2+z^{s-1}\left[1+\left(u_0+v_0+2\right)z-z^2\right] \nonumber
\end{eqnarray}
La famille des polynômes
\[ Q(z)-zP(z) \]
Représente pour $s\geq s_0$, $s_0$ convenablement choisi, une suite d’éléments de l’ensemble $T$ convergente vers le zéro du polynôme
\[1+\left(u_0+v_0+2\right)z-z^2\]
et appartenant à l’ensemble $S$.\newline
En conséquence, pour $s=h$, seul le plus petit élément de $S^\prime(u_0,v_0)$ vérifie que le polynôme $ Q(z)+\eta z^m A$ admet tous ses zéros sur $\left|z\right|=1$.

\subsubsection{Impact de $S^\prime(1)$ sur les points d’accumulation de l’ensemble $T$}

\paragraph{1/ Sur les points d’accumulation appartenant à l’ensemble $S$}\hfill\\
De la suite ($\alpha_s$) de $S^\prime(1)$ représentée par la famille de fonctions de rang infini où
\[A_1(z)=1-z+z^{s-1}-z^s\ ;Q(z)=1-2z+z^{s-1}-z^s \]
Formons les polynômes 
\begin{eqnarray*}
Q(z)+z A_1(z)=1-z-z^2+z^{s-1}\left(1-z^2\right),\ s\geq2\\
Q(z)+z^3 A_1(z)=1-2z+z^2-z^4\left(1-z+z^2\right),\ s=3
\end{eqnarray*}
De la suite ($\beta_s$) de $S^\prime(1)$ représenté par la famille de fonctions de rang infini où
\[\left(1-z\right)A_1(z)=1-z^{s-1}\ ;\left(1-z\right)Q(z)=1-2z+z^{s+1}\]
Formons le polynôme
\[Q(z)+z^3A_1(z)=1-z-z^2+z^{s+1},\ s\geq3\]
De la composition des trois derniers polynômes on établit de (\ref{equ25}) que le point $\tau^\prime$, d’accumulation de l’ensemble $T$, est le zéro de l’un des polynômes
\begin{align}
\label{equ29}
1&-z^2+z^{s-1}\left(1+z-z^2\right),\ s\geq2\nonumber\\
1&-z^{s-1}\left(1+z-z^2\right),\ s\geq3\\
1&-z+z^2-z^4(1-2z+z^2)\nonumber
\end{align}
De Dufresnoy-Pisot \cite{Dufresnoy1955} seuls les zéros des polynômes (\ref{equ29}) appartiennent à l’ensemble $S$ tout en étant inférieurs à l’élément $\theta^\prime$ :
\emph{« En conséquence tout point d’accumulation de l’ensemble T, appartenant à l’ensemble S, est soit le zéro d’un des polynômes (\ref{equ29}), soit supérieur à l’élément $\theta^\prime$. »}\newline

\paragraph{2/ Sur le point $\tau^\prime=1$}\hfill\\
De la suite ($\beta_s$) appartenant à l'ensemble $S^\prime (1)$ et représentée par les polynômes 
\[\left(1-z\right)A(z)=1-z^{s-1}\ ;\left(1-z\right)Q(z)=1-2z+z^{s+1}\]
on obtient
\begin{align}
\label{equ30bis}
Q(z)+zA(z)&=P(z)-zB(z)=1-z^s\nonumber\\
&\qquad\qquad\qquad\qquad \qquad \qquad s\geq2\\
Q(z)+zP(z)&=1-z^{s+1}\nonumber
\end{align}
Le point $\tau^\prime=1$, zéro commun des polynômes (\ref{equ30bis}), ne peut être d’accumulation pour l’ensemble $T$. En conséquence, tout point d’accumulation de l’ensemble $T$, associé à l’ensemble $S^\prime(1)$, appartient à l’ensemble $S$.

\subsubsection{Impact de $S^\prime(u_0)$ sur les points d’accumulation de l’ensemble $T$}\hfill\\
1/ Du théorème \ref{th1} de la représentation des suites ($\gamma_s$) et ($\alpha_s$) par les différents polynômes 
\[Q(z)+z P(z)\]
desquels on forme deux suites d'éléments de l'ensemble $T$, convergentes chacune vers un élément de l'ensemble $S$, à savoir les nombres $\theta_0$ et $\theta^\prime$. \newline
2/ Du théorème \ref{th1} et de la représentation de la suite ($\beta_s$) par les polynômes
\begin{align*}
Q(z)+zA(z)&=P(z)-B(z)=1-z^s \\
Q(z)+zP(z)&=1-z^{s+1}
\end{align*}
l'élément $\tau^\prime=1$ ne peut être d’accumulation de l’ensemble $T$. En conséquence :\newline
\emph{Tout point d’accumulation de l’ensemble $T$, associé à l’ensemble $S^\prime(u_0)$ appartient à l’ensemble $S$.}

\subsubsection{Impact de $S^\prime(u_0,v_0)$ sur les points d’accumulation de l’ensemble $T$}
\paragraph{1/ Pour $s=h$ et $m=1$, $\eta=1$}\hfill\\
Seul le plus petit élément de $S^\prime(u_0,v_0)$, de degré $s=h=2$, de fonction de rang infini $A(z)/Q(z)$, vérifiant
\begin{align*}
A(z)&=u_0-\left(v_0\right)z+z^2\ ;Q(z)=1-\left(u_0+1\right)z+v_0z^2\\
Q(z)&+zA(z)=B(z)+zP(z)=\left(1-z\right)^2(1+z)\\
G(z)&=(u_0-\nu_0)(1-z^2)
\end{align*}
appartient à l’ensemble $S^\prime(u_0,v_0)$.
\paragraph{2/ Sur le point $\tau^\prime$=1}\hfill\\
Du même procédé que dans la partie \ref{I-2}, le point $\tau^\prime=1$ est un zéro commun des polynômes
\begin{align}
\label{equ31}
Q(z)+zA(z)&=P(z)-B(z)=1-z^s \nonumber\\
& \qquad  \qquad  \qquad  \qquad  \qquad   \qquad  s\geq2\\
Q(z)+zP(z)&=1-z^{s+1}\nonumber
\end{align}
formés de la représentation de la suite ($\beta_s$) du théorème \ref{th1}, $s=h+1$.\newline
Le point $\tau^\prime=1$, zéro commun des polynômes (\ref{equ31}), ne peut être d’accumulation de l’ensemble $T$. En conséquence :\newline
\emph{Tout point d’accumulation de l’ensemble $T$, associé à l’ensemble $S^\prime(u_0,v_0)$ appartient à l’ensemble $S$.}

De cette étude on établit
\begin{lemme} \label{lemme7}
Tout point d’accumulation de l’ensemble $T$ appartient à l’ensemble $S$.
\end{lemme}

\section{Conclusion}
Du Lemme \ref{lemme7} associé au Lemme 5 de Salem dans \cite{Salem1944}, on conclut notre étude par le résultat :\newline
L’ensemble $S\cup T$ est un fermé de la demi-droite réelle $\left]1, +\infty\right[$ de plus petit élément $\tau_0$, zéro du polynôme \[1+z-z^3\left(1+z+z^2+z^3+z^4\right)+z^9+z^{10}\]
rejetant ainsi la conjecture de Lehmer \cite{Lehmer1933} et confirmant celle de Boyd \cite{Boyd1977}.\newline

Je souhaite dédier ce travail à la mémoire de mon maître Charles Pisot et tiens à remercier mon collègue Mongi Naïmi, pour l'intérêt manifesté tout au long de cette étude et son aide dans l'élaboration du manuscrit.


\bibliographystyle{abbrv}

\bibliography{/Applications/TeX/Lehmer}

%
%
%
%
%

\end{document}